\documentclass[12pt]{article}
\usepackage{amsmath,verbatim,amsfonts,mathrsfs}
\usepackage{verbatim}% for comment
\def\bbbr{{\rm I\!R}}

\newcommand{\Proof}{\noindent{\it Proof. }}
\newcommand{\Endproof}{\Box\medskip}

\newcommand{\bi}{\mathbf{i}}
\newcommand{\bj}{\mathbf{j}}

\newcommand{\scrC}{\mathscr{C}}
\newcommand{\dimh}{\dim_{\rm H}}
\newcommand{\dimp}{\dim_{\rm P}}

\renewcommand{\Box}{\mbox{\rule{1ex}{1ex}}}
\newtheorem{pro}{Proposition}[section]
\newtheorem{lem}[pro]{Lemma}
\newtheorem{thm}[pro]{Theorem}

\newcounter{example}
\newcommand\hdd{\mbox{\rm dim}_{\rm H}\,} % Hausdorff dimension
\newcommand\ubd{\overline{\mbox{\rm dim}}_{\rm B}\,} % upper box dimension
\newcommand\lbd{\underline{\mbox{\rm dim}}_{\rm B}\,} % lower box dimension
\newcommand{\ssi}{_{\bf i}} % sub bf i
\newcommand{\w}{\omega} % random translate
\newcommand{\E}{{\sf E}} % expectation
\renewcommand{\P}{{\sf P}} % probabiltiy
\newcommand{\be}{\begin{equation}} % begin equation
\newcommand{\ee}{\end{equation}} % end equation

\title{Local dimensions of measures on self-affine sets}
\author{K.J. Falconer and Jun Jie Miao \footnote{This work is supported by the
National Natural Science Foundation of China 10971069 and Shanghai
Education Committee Project 11ZZ41}\\}
\begin{document}
\maketitle

\begin{abstract}
We show that, in a generic setting, self-affine and almost self-affine measures are exact dimensional, with local dimension equal almost everywhere to the information dimension and given by the zero of a superadditive pressure functional. 
\end{abstract}

\section{Introduction}
\setcounter{equation}{0}
\setcounter{pro}{0}

For $m \geq 2$ let $T_{1},\ldots,T_{m}$ be a set of linear
contractions on $\bbbr^{N}$, let $\w_{1},\ldots,\w_{m} \in \bbbr^{N}$
be a set of translation vectors, and let
$S_{1},\ldots,S_{m}:\bbbr^{N}\to \bbbr^{N}$ be the  affine transformations
$$
    S_{i}(x) = T_{i}(x) +\w_{i} \qquad (i=1,\ldots,m).
$$
The contractions $\{S_{1},\ldots,S_{m}\}$ form an {\em
iterated function system} (IFS).  By the well-known theorem of
Hutchinson, see \cite{Fa,Hut}, this IFS has a unique {\em
attractor}, that is a unique non-empty compact set  $E^\w \subset
\bbbr^{N}$ such that
$$
E^\w = \bigcup_{i=1}^{m} S_{i}(E^\w).
$$
 We refer to  $E^\w$
as  a {\em self-affine} set, writing $E^\w$  to emphasise its dependence on the
vector of translations $\w=(\w_1,\ldots,\w_m) \in \bbbr^{Nm}$.

Whilst a great deal is known about the dimensions of self-similar sets, see for example \cite{Fa,Hut}, determination of the Hausdorff and box dimensions of self-affine
sets can be challenging, not least because the dimensions need not vary continuously with $(\w_1,\ldots,\w_m)$.
A covering argument, involving dividing up
the components of $E^\w$ into appropriate pieces, shows that, for all
$\w$, the Hausdorff and lower and upper box-counting dimensions
satisfy
\begin{equation}\label{ineq_dim}
\hdd (E^\w ) \leq \lbd (E^\w) \leq \ubd (E^\w) \leq
\min\{d_0,N\} ,
\end{equation}
where $d_0$, termed the {\it affinity dimension}, satisfies 
$$
 \lim_{k\to \infty}\frac{1}{k}
\log \sum_{|\bi|=k} \phi^{d_0}(T_\bi)= 0
$$
(see Section  4 for the notation used here).

%\begin{equation}\label{dimexp}
%d(T_{1},\ldots, T_{m}) = \inf \Big\{s: \sum_{k=1}^{\infty} \sum_{\bi
%\in \bJ_{k}}\svs {T_{\bi}} < \infty\Big \};
%\end{equation}
%see\cite{Fa1,JPS}.

For many self-affine sets  equality holds in (\ref{ineq_dim}), at least in a `generic' sense, see  \cite{Edg,Fa1, Fa2,HL,Kae,KS, PSol,Sol}.
In particular this is true for ${\cal L}^{mN}$-almost all $(\w_1,\ldots,\w_m)
\in \bbbr^{mN}$, provided $\|T_i\| < \frac{1}{2}$ for all $i$, see
\cite{Fa1, Sol}. Nevertheless,   for many regular
self-affine sets such as affine Sierpi\'{n}ski  carpets (where the $S_{i}$ map the unit
square onto rectangles selected from a rectangular grid) the
dimension is in general strictly less than
 the affinity dimension, see \cite{Bar,Bed,FL,FW,LG,McM}. 

A variant of self-affine sets  with rather more randomness was recently introduced \cite{JPS}. Here there is a scaled 
independent random translation at each stage of the
iterated construction of the set to yield a statistically
self-affine set, and this was shown to have Hausdorff and box
dimensions $\min\{d_0,N\} $ almost surely, provided only that $\| T_{i} \| < 1$ for all $i$.
% (Note that it is often convenient to use the
%language of probability rather than of measure theory when
%considering such constructions.)

It is natural to consider multifractal aspects of measures supported by these sets. In particular, self-affine measures, which may be thought of as Bernoulli measures on a code space projected onto self-affine sets in the natural way, have been studied in various cases. Again, multifractal quantities need not vary continuously with the defining parameters, and there are two approaches. One can consider measures on specific self-affine sets, such as affine Sierpi\'{n}ski carpets, see  \cite{ BM, K,O}. Alternatively, one can seek generic formulae valid almost surely across some parameter space. In particular,  the almost sure generalized $q$-dimensions of measures on self-affine sets  \cite{Fa5}  and  on almost self-affine sets \cite{Fa6}, have been obtained  for certain ranges of $q$.

In this paper we consider the pointwise or local dimension of measures, that is
$$\lim_{r\to 0}\ \frac{\log \mu^\w (B(x,r))}{\log r}.$$
We show  that in a generic situation $\mu^\w$ is exact dimensional, that is, for almost all $\w$ in some parameter space,  the local dimension of $\mu^\w$ exists and takes a common value $d_1$ at $\mu^\w$-almost all $x$, where  $d_1$ is the zero of a superadditive pressure functional (\ref{dqdef}) and equals  the information dimension  ${D}^1(\mu^\w)$.

To establish this, we obtain a deterministic upper bound for the local dimensions by utilizing the upper bounds for the $L^q$-dimensions as $q\nearrow 1$ together with the lower semicontinuity of a family of pressure functions $P(s,q)$ at $q=1$, see Proposition \ref{thm_ub}. The almost sure lower bound comes from an integral estimate, Theorem \ref{thm_lb}.

We end with estimates in the more general case when $\mu$ is a Gibbs measure on the code space.

\section{Dimensions of measures}
\setcounter{equation}{0}
\setcounter{pro}{0}

In this section we review the multifractal definitions and properties that we will require; see \cite{Fa3,P} for full accounts.

We denote the closed ball of radius $r$ with center $x$ by $B(x,r)$.
Let $\nu$  be a Borel regular probability measure on $\bbbr^N$. The
\textit{Hausdorff dimension} and \textit{packing dimension} of $\nu$
are defined by
\begin{equation*}
\dimh \nu=\sup\Big\{s:\liminf_{r\to 0}\frac{\log \nu (B(x,r))}{\log
r}\geq s, \qquad\textrm{for $\nu$-almost all } x\Big\},
\end{equation*}
\begin{equation*}
\dimp \nu=\sup\Big\{s:\limsup_{r\to 0}\ \frac{\log \nu (B(x,r))}{\log
r}\geq s, \qquad\textrm{for $\nu$-almost all } x\Big\}.
\end{equation*}
If for some $s$
$$\lim_{r\to 0}\ \frac{\log \nu (B(x,r))}{\log r} = s,$$
(with the limit existing) for $\nu$-almost all $x$,
we say that $\nu$ is {\it exact of dimension} $s$, in which case $s=\dimh \nu=\dimp \nu$.

We write $\mathcal{M}_r$ for the family 
of $r$-\textit{mesh cubes} in $\bbbr^N$, that is cubes of the form
\newline
$[j_1 r, (j_1+1)r)\times \cdots\times [j_N r, (j_N+1)r)$ where $j_1,
\ldots, j_N \in \mathbb{Z}$. Let $\nu$  be a Borel probability measure on $\mathbb{R}^N$.
For $q \neq 1$ we define the {\it lower} and {\it upper generalized $q$-dimensions} or $L^q$-{\it dimensions} of $\nu$ to be
\be
\underline{D}^q(\nu)=\liminf_{r\to 0}\frac{\log \sum_{\mathcal{M}_r}\nu(C)^q}{(q-1)\log r}, \quad
\overline{D}^q(\nu)=\limsup_{r\to 0}\frac{\log \sum_{\mathcal{M}_r}\nu(C)^q}{(q-1)\log r}.
\label {gendim}
\ee
For $q=1$, $\underline{D}^1(\nu)$ and $\overline{D}^1(\nu)$, also termed the
\textit{lower} and \textit{upper information dimensions}, are defined by
\be\label {infdim}
\underline{D}^1(\nu)=\liminf_{r\to 0}\frac{\sum_{\mathcal{M}_r}\nu(C)\log\nu(C)}{\log r}, \quad
\overline{D}^1(\nu)=\limsup_{r\to 0}\frac{\sum_{\mathcal{M}_r}\nu(C)\log\nu(C)}{\log r}.
\ee
If $\underline{D}^q(\nu)=\overline{D}^q(\nu)$,
we write $D^q(\nu)$ for the common value which we refer to as the
\textit{generalized $q$-dimension}.

We will need the following lemma  on several occasions.

\begin{lem}\label{lem_cts}
Let $\nu$  be a probability measure on some measure space $X$ and let $f: X \to \bbbr^+$ be $\nu$-measurable. Write
\begin{eqnarray*}
 F(q) &=& \frac{1}{q-1}\log \int f^{q-1} d \nu  \quad (q \neq 1)\\
 F(1) &=&  \int \log f d \nu  
 \end{eqnarray*}
 Then $F: \bbbr \to \bbbr\cup\{-\infty,\infty\}$ is a monotonic increasing function.
 Moreover, if $0<a \leq f(x) \leq b<\infty $ for all $x \in X$ for some $a$ and $b$, then $F: \bbbr \to \bbbr$ is continuous.
\end{lem}
\Proof Monotonicity of $F$ follows from Jensen's inequality.
If $f$ is bounded away from $0$ and $\infty$, continuity at $q \neq 1$ follows from the bounded convergence theorem. 
For $q$ close to $1$, note that
\begin{eqnarray*}
 \int f^{q-1} d \nu &=& \int\exp\big((q-1)\log f\big) d\nu =  \int \big(1+ (q-1) \log f +O((q-1)^2)\big) d \nu \\
 &=& 1 + (q-1) \int \log f  d\nu +O((q-1)^2),
\end{eqnarray*}
so 
$$\log \int f^{q-1} d \nu = (q-1) \int \log f  d\nu +O((q-1)^2),$$
giving continuity of $F(q)$ at $q=1$.
\Endproof

It is often convenient to express generalized dimensions as integrals of measures of balls rather than as moment sums.

\begin{pro}
The generalized dimensions have integral forms:
\be\label{gendimint}
\underline{D}^q(\nu)=\liminf_{r\to 0}\frac{\log\int\nu(B(x,r))^{q-1} d\nu(x)}{(q-1)\log r}, \quad
\overline{D}^q(\nu)=\limsup_{r\to 0}\frac{\log\int\nu(B(x,r))^{q-1} d\nu(x)}{(q-1)\log r},
\ee
for $q>0, q \neq 1$, and
\be\label{infdimint}
\underline{D}^1(\nu)=\liminf_{r\to 0}\frac{\int\log\nu(B(x,r)) d\nu(x)}{\log r}, \quad
\overline{D}^1(\nu)=\limsup_{r\to 0}\frac{\int\log\nu(B(x,r)) d\nu(x)}{\log r}.
\ee
Moreover, $\underline{D}^q(\nu)$ and $\overline{D}^q(\nu)$ are monotonic decreasing in $q$.

\end{pro}

\Proof
Identity (\ref{gendimint}) is straightforward for $q>1$, see, for example, \cite{P}. The case of $0<q<1$ was established in \cite{PSol1}.

For $q=1$, given $0<r<1$, for each $x\in \bbbr^N$ write $C(x)$ for the $r$-mesh cube containing $x$. Then
$$
\sum_{\mathcal{M}_r}\nu(C)\log\nu(C) 
= \int\log\nu(C(x)) d\nu(x)
\leq \int\log\nu(B(x,\sqrt{N} r) d\nu(x),
$$
and dividing by $\log r$ and taking the limits gives that the expressions of (\ref{infdim}) are  at least the corresponding ones of  (\ref{infdimint}).

For the opposite inequalities, fix  $0<r<1$ and write $\widetilde{C}$ for the cube of side $3r$ formed by the $3^N$ cubes in $\mathcal{M}_r$ consisting of $C$ and its immediate neighbours. 
Let ${\cal S}_k \; (k=1,2,3,\ldots)$ be the set of mesh cubes
$${\cal S}_k = \{C \in \mathcal{M}_r: 2^{k-1}\nu(C) \leq \nu(\widetilde{C} )< 2^{k}\nu(C)\}.
$$
Then
\be
\sum_{C \in \mathcal{S}_k}\nu(C)
\leq 2^{1-k} \sum_{C \in \mathcal{S}_k}\nu(\widetilde{C})
\leq 2^{1-k} 3^N \sum_{C \in  \mathcal{M}_r}\nu(C) = 2^{1-k} 3^N. \label{cksize}
\ee
Writing $\widetilde{C}(x)$ for the cube $\widetilde{C}$ containing 
$x$,
\begin{eqnarray*}
\int\log\nu(B(x,r)) d\nu(x) 
&\leq& \int\log\nu(\widetilde{C}(x)) d\nu(x) \\
&=& \sum_{C \in \mathcal{M}_r}\nu(C) \log\nu(\widetilde{C})\\
&\leq& \sum_{k=1}^\infty \sum_{C \in \mathcal{C}_k}\nu(C) \log(2^k\nu(C))\\
&\leq& \sum_{k=1}^\infty \sum_{C \in \mathcal{C}_k}\nu(C)\big(\log\nu(C) +k\log 2\big)\\
&=& \sum_{C \in \mathcal{M}_r}\nu(C)\log \nu(C) +   \sum_{k=1}^\infty 3^N 2^{1-k}k \log 2,
\end{eqnarray*}
using (\ref{cksize}). Since the right hand sum is finite, dividing by $\log r$ and taking the limit completes the argument for $q=1$.

Monotonocity of $\underline{D}^q(\nu)$ and $\overline{D}^q(\nu)$ follows from Lemma \ref{lem_cts}.
\Endproof

The local dimensions of a measure are related to the generalized dimensions as $q$ approaches $1$.

\begin{pro}\label{propat1}
Let $\nu$  be a Borel probability measure on  $\bbbr^N$. Then for $\nu$ almost-all $x$
\be
\lim_{q \searrow 1} \underline{D}^q(\nu)
\leq \liminf_{r\to 0} \frac{\log \nu(B(x,r))}{\log r}
\leq \limsup_{r\to 0} \frac{\log \nu(B(x,r))}{\log r}
\leq \lim_{q \nearrow 1} \overline{D}^q(\nu).\label{mesdimineq}
\ee
\end{pro}

\Proof
For the right hand inequality, let $q<1$ and let  $ \overline{D}^q(\nu)<t<s$. From (\ref{gendimint}) there is a constant $c$ such that 
$$\int \nu(B(x,r))^{q-1} d\nu(x) \leq c r^{(q-1)t}$$
for all $0<r\leq 1$,
so that 
$$\int \left( \frac{\nu(B(x,r))}{r^s}\right)^{q-1} d\nu(x) \leq c r^{(q-1)(t-s)}.$$
Setting $r= 2^{-k}$ for $k=1,2,\ldots$ and summing, we get
$$\int \sum_{k=1}^{\infty}
 \left(\frac{\nu(B(x,2^{-k}))}{2^{-ks}}\right)^{q-1} d\nu(x) \leq c  \sum_{k=1}^{\infty}2^{-k(1-q)(s-t)}<\infty.$$
It follows that, for $\nu$-almost all $x$, 
$\lim_{k\to \infty} \nu(B(x,2^{-k}))/2^{-ks} = \infty$. By comparing $\nu(B(x,r))$ with $\nu(B(x,2^{-k}))$ where $2^{-k} <r \leq 2^{-k+1}$, this implies that  
$\lim_{r\to 0} \nu(B(x,r))/r^s = \infty$ so 
$\limsup_{r\to 0} \log \nu(B(x,r))/\log r \leq s$. 
This is so for all $s> \overline{D}^q(\nu)$, giving  the right hand inequality.

The left hand inequality is similar.
\Endproof

%333333333333333333333333

\section{Almost self-affine sets and measures}
\setcounter{equation}{0}
\setcounter{pro}{0}

In this section we recall the code space representation of almost self-affine sets, of  which self-affine sets are a special case.

As is usual, we index a subset of  $\bbbr^N$ constructed in a hierarchical manner by
a {\it code space} or {\it sequence space}. Then a measure on the code space may be projected to a measure on the set. Let $m \geq 2$.
For $k=0,1,2, \ldots$  let $I_{k}$ be the
set of all $k$-term sequences of integers $1,2, \ldots, m$,
that is $I_{k} = \{ (i_{1}, \ldots , i_{k}): \, 1 \leq i_{j} \leq 
m \}$;
we regard $I_{0}$ as just containing the empty sequence $\emptyset$. 
 We abbreviate members of $I_{k}$ by
${\bf i} = (i_{1},  \ldots , i_{k} ) $ and write $|{\bf i}| = k$ for 
the number of terms in ${\bf i}$.
We write
$I= \cup^{\infty}_{k=0} I_{k} $
for the set of all such finite sequences, and $I_{\infty}$ for the corresponding
set of infinite sequences, so
$I_{\infty} = \{(i_{1}, i_{2}, \ldots ): 1 \leq i_{j} \leq m \}$.
Juxtaposition
of ${\bf i}$ and ${\bf j}$ is denoted by ${\bf ij}$.  
We write ${\bf 
i}|k = (i_{1}, \ldots , i_{k})$ for the {\it curtailment} after $k$ terms
of ${\bf i} = (i_{1}, i_{2}, \ldots ) \in I_{\infty}$, 
with a similar notation if ${\bf i} \in I_{k^{\prime}}$ and $ k 
\leq k^{\prime}$.  We write ${\bf i} \preceq{\bf j}$ if ${\bf i}$ is a 
curtailment of ${\bf j}$.  If ${\bf i,j} \in I_{\infty} $  then 
${\bf i}\wedge {\bf j}$ denotes the common initial subsequence of ${\bf i} $ and ${\bf j}$, that is the maximal sequence such that both 
${\bf i}\wedge {\bf j} \preceq{\bf i}$ and ${\bf i}\wedge {\bf j} \preceq {\bf j}$.

We topologise $I_{\infty}$  using the metric 
$d(\bi,\bj) = 2^{-|\bi \wedge \bj |}$ for distinct $\bi,\bj \in 
I_{\infty}$ to make $I_{\infty}$  a compact metric space. 
The {\it cylinders}  
$C \ssi = \{\bj \in I_{\infty} : \bi \preceq \bj \}$ for  $\bi \in I$ form a base 
of open and closed neighbourhoods for $I_{\infty}$.

Let $T_{1},\ldots,T_{m}$ be a set of linear
contractions on $\bbbr^{N}$. For each $\bi = (i_1,\ldots,i_k) \in I_k$ let   $\w_{\bi} = \w_ {i_1,\ldots,i_k} \in \bbbr^{N}$
be a  translation vector, and let $\w = \{\w_{\bi}: \bi \in I\}$ denote the family of such translations.
We assume throughout that there is some non-empty compact set $B \subset \bbbr^N$ such that 
\be
T_{i_1}(B) + \w_ {i_1,\ldots,i_k} \subseteq B
\label{inclusions}
\ee
for all  $\bi = (i_1,\ldots,i_k) \in I_k$ and for all $\w$ under consideration. This ensures that each 
$\bi = (i_1,i_2, \ldots) \in I_\infty$ determines a nested set of affine copies of $B$ with intersection the single point
\begin{eqnarray}
x^\w(\bi) & =&\bigcap_{k=0}^\infty (T_{i_1} + \w_{i_1})  (T_{i_2} +
\w_{i_1,i_2}) \cdots
  (T_{i_k} + \w_{i_1,\ldots,i_k})(B)\label{points1}\\
  & =&\lim_{k \to \infty} (T_{i_1} + \w_{i_1})  (T_{i_2} +
\w_{i_1,i_2}) \cdots
  (T_{i_k} + \w_{i_1,\ldots,i_k})(x)\label{points2}\\
& =&  \w_{i_1} +T_{i_1}\w_{i_1,i_2} + T_{i_1}T_{i_2}\w_{i_1,i_2,i_3}
+ \cdots. \label{points3}
\end{eqnarray}
It is easily checked that, given (\ref{inclusions}), these limits exist, that   (\ref{points2})  is independent of $x\in \bbbr^N$, and that the map $\bi \mapsto x^\w(\bi)$ is continuous for every $\w$.

We term the compact set $E^{ \w}$ given by the aggregate of these points, 
\begin{equation}
E^{ \w} = \bigcup_{{\bf i} \in I_{\infty}} x^\w(\bi) \subset
\bbbr^N, \label{2.6}
\end{equation}
an {\it almost self-affine set}, and this may be written as
\begin{equation}
E^{ \w} = \bigcap_{k=0}^{\infty} \bigcup_{i_1,\dots,i_k \in I_k}
(T_{i_1} + \w_{i_1}) (T_{i_2} + \w_{i_1,i_2}) \cdots (T_{i_k} +
\w_{i_1,\ldots,i_k})(B), \label{construct}
\end{equation}
which represents the standard hierachical way of constructing a fractal set $E^\w$, see \cite{Fa6,JPS}.

In the special case where 
$ \w_ {i_1,\ldots,i_k} =  \w_{i_k}$ for all $\bi = (i_1,\ldots,i_k) \in I_k$, the set $E^{\w}$ is the unique non-empty compact subset of $\bbbr^N$ satisfying
\begin{equation}
E^{\w} = \bigcup_{i=1}^{m} S_{i}(E^{\w})\label{attractor}
\end{equation}
where $S_{i}(x) = T_{i}(x) +\w_{i}\,  (i=1,\ldots,m)$ are contracting  affine  transformations,  so $E^{\w}$ is the attractor of the iterated function system  $\{S_1,\ldots, S_m\}$ and    $E^{\w}$
is  a {\em self-affine} set. (By a slight abuse of notation we write $\w = (\w_1,\ldots,\w_m)$ for the translation parameters in this situation.)

Let $\mu$  be a finite Borel regular measure on $I_\infty$. For each $\w$, we define  $\mu^\w$, the projection of $\mu$ measure onto $\bbbr^N$, by
\begin{equation}
 \mu^\w(A)=\mu\{\bi: x^\w(\bi)\in A\}                    \label{def_ua}
\end{equation}
for $A\subseteq\bbbr^N$, or equivalently by
\be
\int f(x)d\mu^\w(x)=\int f(x^\w(\bi))d\mu(\bi) \label{def_ua2}
\ee
for every continuous  $f:\bbbr^N\to \bbbr$. Then $\mu^\w$ is a Borel measure supported by $E^{\w}$.

In particular, given `probabilities' $p_1,\ldots,p_m$ (so that $p_i>0$ for each $i$ and $\sum_{i=1}^m p_i=1$) we may define a Bernoulli measure $\mu$ on $I_\infty$ by setting
\begin{equation}\label{mP}
\mu(C_\bi)=p_\bi\equiv p_{i_1}p_{i_2}\ldots p_{i_k} \quad (\bi=i_1\ldots i_k)
\end{equation}
for each cylinder $C_\bi$ and
extending to general subsets of $I_\infty$ in the usual way.  With $E^\w$ an almost self affine set, the projected measure $\mu^\w$ on $E^{
\w}$  given by (\ref{def_ua}) and (\ref{def_ua2})  is termed an {\it almost self-affine measure}, and if $E^{\w}$ is a self-affine set (\ref{attractor})  $\mu^\w$ is termed a {\it self-affine measure}, in which case 
\begin{equation}\label{msa}
\mu^{\w}(A)=\sum_{i=1}^m p_i \mu^{\w} (S_i^{-1}(A)),
\end{equation}
for $A \subseteq \bbbr^N$, see \cite{Fa3,Hut}.
 
%444444444444444444444444

\section{The pressure functions}
\setcounter{equation}{0}
\setcounter{pro}{0}

The fractal and multifractal behaviour of self-affine and almost self-affine sets and measures depend on certain `subadditive pressure functionals' defined in terms of singular value functions of the linear mappings $T_i$, see \cite{Barr,Fa1,Fa5,Fa6, Sol}.
 
Let $T:\bbbr^N\rightarrow\bbbr^N$ be a non-singular linear
contraction. The {\em
singular values} $\alpha_i \equiv \alpha_i(T)$ of $T$
($i=1,\ldots,N$) are the positive square roots of the eigenvalues of
$T^{*}T$, where $T^{*}$ is the transpose or adjoint of $T$. Equivalently they
are the lengths of the principal semi-axes of the image $T(B)$ of
the unit ball $B$. We adopt the convention that
$1>\alpha_1\geq\alpha_2\geq\cdots\geq\alpha_N>0$. The {\em singular
value function} $\phi^s(T)$ is then defined for $ 0 \leq s \leq N$
as
$$
\phi^s(T)=\alpha_1\alpha_2\cdots\alpha_{m-1}\alpha_m^{s-m+1},
$$
where $m$ is the integer such that $m-1< s\leq m$, with the
convention that $\phi^s(T)=(\alpha_1\alpha_2\cdots\alpha_{N})^{s/N}$
if $s \geq N$. The singular value function $\phi^s(T)$ is decreasing
in $s$ and is submultiplicative, that is
$\phi^s(T_{\bi\bj})\leq\phi^s(T_{\bi})\phi^s(T_{\bj})$ for all
$\bi,\bj \in I$, where we write $T_{\bi} = T_{i_1}T_{i_2}\cdots T_{i_k}$ for
$\bi= i_1,\ldots,i_k \in I$.

Let
$$
\alpha_{+} = \max_{i=1,\ldots,N}\alpha_{1}(T_{i}) \quad \mbox{ and
}\quad \alpha_{-} = \min_{i=1,\ldots,N}\alpha_{n}(T_{i}).
$$
This gives the bounds
$$
\alpha_{-}^{s|\bi|} \leq \phi^{s}(T_{\bi})\leq \alpha_{+}^{s|\bi|}
\qquad ( \bi \in I).
$$

The following expressions, which may be thought of as subadditive pressure functionals, are central to the theory of self-affine measures. It is helpful to view them both as sums over cylinders and as integrals.

For $s \geq 0$ and $q \geq 0, q \neq 1$ define:
\begin{equation}
P(s,q) = \lim_{k\to \infty}\frac{1}{k}
\frac{\log \sum_{|\bi|=k} \phi^s(T_\bi)^{1-q}\mu(C_\bi)^q}{q-1}
\equiv \lim_{k\to \infty}\frac{1}{k}\frac{\log \int \big(\phi^s(T_{\bi|k})^{-1}\mu(C_{\bi|k})\big)^{q-1}d\mu(\bi)}{q-1},
\label{Pdef}
\end{equation}
and for $s \geq 0$ and $q = 1$,
\begin{equation}
P(s,1) = \lim_{k\to \infty}\frac{1}{k}
 \sum_{|\bi|=k} \mu(C_\bi)\log \big( \phi^s(T_\bi)^{-1}\mu(C_\bi)\big) 
\equiv \lim_{k\to \infty}\frac{1}{k}
\int \log \big( \phi^s(T_\bi)^{-1}\mu(C_{\bi|k})\big) d\mu(\bi). \label{P1def}
\end{equation}
These limits exist for a Bernoulli measure $\mu$ since, from the submultiplicativity of the $\phi^s(T_\bi)$, the sums in (\ref{Pdef}) form a sub- or supermultiplicative sequence in $k$ (depending on whether $0<q<1$ or $q>1$), and the sum in (\ref{P1def}) is superadditive. We note the following properties of $P(s,q)$.

\begin{lem}\label{lemPcts}
For each $q\geq0$, $P(s,q)$ is strictly monotonic increasing and continuous in  $s$. More precisely,
\be
0 < h \log \alpha_+^{-1} \leq P(s+h,q) -P(s,q) \leq h \log \alpha_{-}^{-1}
\quad (s \geq 0, h>0). \label{scts}
\ee

For each $s\geq 0$, $P(s,q)$ is monotonic increasing in $q$, is lower semicontinuous in $q$ and continuous for  $q\neq 1$. Moreover, there are $0<\beta_{-}\leq \beta_{+}<\infty$ such that
\be
(q'-q) \beta_{-} \leq (q'-1)P(s,q') -(q-1)P(s,q) \leq (q'-q) \beta_{+}
\quad (0\leq q \leq q'<1 \mbox{ \rm or } 1< q \leq q'). \label{pcts}
\ee
\end{lem}

\Proof
Inequalities (\ref{scts}) follow from the definitions (\ref{Pdef}),(\ref{P1def}) of
$P(s,q)$, noting that 
$\phi^s(T_\bi) \alpha_{-}^{kh} \leq \phi^{s+h}(T_\bi) \leq \phi^s(T_\bi) \alpha_{+}^{kh}$.

Fixing $s \geq 0$,  Lemma \ref{lem_cts} gives that for each $k$, 
$\log \sum_{|\bi|=k} \phi^s(T_\bi)^{1-q}\mu(\scrC_\bi)^q/(q-1) = \log \int \big(\phi^s(T_{\bi|k})^{-1}\mu(C_{\bi|k})\big)^{q-1}d\mu(\bi)/(q-1)$ (with the logarithmic definition when $q=1$) is continuous and increasing with $q$, so the same is true for the limit  $P(s,q)$. Moreover, these expressions are
superadditive for each $q$, so by a standard property of superadditive sequences, $P(s,q)$ is not only the limit as $k \to \infty$ in  (\ref{Pdef}) and (\ref{P1def}) but also the supremum over $k$. Since these expressions are continuous, $P(s,q)$ is lower semicontinuous as the supremum of a family of continuous functions.

Inequalities (\ref{pcts}), which imply  continuity for $q \neq 1$, follow from the definitions (\ref{Pdef}),(\ref{P1def}) of
$P(s,q)$, taking 
$\beta_{-} = \log(\alpha_{+}^{-s}\min_i \mu(C_i))$ and    
$\beta_{+} = \log(\alpha_{-}^{-s}\max_i \mu(C_i))$.

\Endproof

Note that in certain cases, for example when $\phi^s(T_\bi)$ is multiplicative, then $P(s,q)$ may be continuous for all $q\geq 0, s\geq 0$; this  happens if the $T_i$ are similarities or if the $T_i$ can all be represented by diagonal matrices with respect to some basis.

The previous lemma guarantees that, for each $q\geq 0$, there is a unique number $d_q>0$ such that 
$P(d_q,q) = 0$; specifically the $d_q$ satisfy
\begin{eqnarray}
P(d_q,q) &=& \lim_{k\to \infty}\frac{1}{k}
\frac{\log \sum_{|\bi|=k} \phi^{d_q}(T_\bi)^{1-q}\mu(C_\bi)^q}{q-1} = 0
\qquad (q \neq 1) \nonumber\\
 P(d_1,1) &=& \lim_{k\to \infty}\frac{1}{k}
\log \sum_{|\bi|=k} \mu(C_\bi)\log \big( \phi^{d_1}(T_\bi)^{-1}\mu(\bi)\big)
= 0 .\label{dqdef}
\end{eqnarray}

\begin{lem}\label{lem_ctsdq}
For  $q\geq0$, $d_q$ is strictly monotonic decreasing in $q$. Furthermore, $d_q$ is continuous at all $q \neq 1$ and is upper semicontinuous at $q=1$.
\end{lem}

\Proof
As $P(s, q)$   is  strictly increasing in $q$ and increasing in $s$, the solution $d_q$ of $P(d_q, q)=0$ is strictly decreasing in $q$.

To show upper semicontinuity of  $d_q$ at $q \geq 0$, let $s> d_{q}$. Then
$P(s, q)>0$, so by lower semicontinuity of $P(s,q)$ in $q$, there exists
$\delta >0$ such that if $|q-q'|<\delta$ then $P(s, q')>0$ and so $s>d_{q'}$ by monotonicity. Thus  $d_q$ is upper semicontinuous at $q$.

A symmetric argument using the upper semicontinuity of $P(s,q)$ in $q$ at $q \neq 1$ shows that $d_q$ is lower semicontinuous at $q \neq 1$.
\Endproof

%5555555555555555555555

\section{Dimensions of measures on almost self-affine sets}
\setcounter{equation}{0}
\setcounter{pro}{0}

We now derive estimates for the local dimensions of
$\mu^\w$. We assume throughout that  the support $E^{\w}$satisfies (\ref{inclusions}) for some non-empty compact $B$.

We recall the following upper bound for generalized $q$-dimensions.

\begin{pro}\label{Dq_ub}
Let $\mu$  be a Bernoulli measure on  $I_\infty$. If $q\geq 0, q \neq 1$ then $\overline{D}^q(\mu^\w) \leq \min\{d_q,N\}$ for all $\w$ such that $(\ref{inclusions})$ is satisfied.
\end{pro}

\Proof
This upper bound is derived in \cite[Proposition 4.1, Theorem 6.2]{Fa5} by subdivision of the sets $T_\bi(B)$ and applying H\"{o}lder's inequality to the parts. \Endproof

Note that Proposition \ref{Dq_ub} is also valid for $q=1$, though we do not need this here.

Taking the limit as $q \nearrow 1$ of the above estimate gives an upper bound for the local dimensions.

\begin{pro}\label{thm_ub}
Let $\mu$  be a Bernoulli measure on $I_\infty$. For any   $\w$ let $\mu^{\w}$ be the projection of $\mu$  onto $E^\w$ given by $(\ref{def_ua})$. Then
$$
\limsup_{r\to 0}\frac{\log\mu^\w(B(x,r))}{\log r}\leq \min\{d_1,N\} 
$$
for $\mu^\w$-almost all $x\in \bbbr^{N}$.
\end{pro}
\Proof   
By Proposition \ref{Dq_ub}  $\overline{D}^q(\mu^\w) \leq d_q$ for all $0<q<1$.
Hence, using the right hand inequality  of (\ref{mesdimineq}), 
$$
\limsup_{r\to 0}\frac{\log\mu^\w(B(x,r))}{\log r}
\leq \lim_{q \nearrow 1}\overline{D}^q(\mu^\w) 
\leq \lim_{q \nearrow 1} d_q
\leq d_1,
$$
using the upper semicontinuity of $d_q$. Recall also that the upper local dimension of any measure on $\bbbr^N$ is at most $N$ almost everywhere.
\Endproof

Given that the dimension of $E^\w$ need not be continuous in the translations $\w$, we can only expect to show that  $\min\{d_1,N\}$ is also a lower bound for local dimensions for almost all constructions, in some sense. It is convenient to express this in probabilistic language. Thus let 
$\Omega = \{\w_\bi : \bi \in I\}$ and let $\P$ be a probability measure on $\Omega$ such that the random vectors $\w_{\bi}$ are jointly measurable. We write $\E$ for expectation with respect to $\P$. Here is a general result which we specialize to more specific probability distributions in the next section.

\begin{thm}\label{thm_lb}
Let $\mu$  be a Bernoulli measure on $I_\infty$ given by $(\ref{mP})$ and
let $\mu^{\w}$ be the projection of $\mu$ onto $E^\w$. 
Suppose that there are numbers $s<\min\{d_1,N\}$ that are arbitrarily close to $\min\{d_1,N\}$ for which
there exists $c<\infty$ such that
\be\label{potest}
\E\left( |x^\w(\bi)-x^\w(\bj)|^{-s}\right)\leq
\frac{c}{\phi^s(T_{\bi\wedge \bj})} \quad (\bi \neq \bj \in I_\infty). 
\ee
Then, for almost all $\w$, the measure $\mu^\w$ is exact dimensional, with
$$
\lim_{r\to 0}\frac{\log \mu^\w (B(x,r))}{\log r} = D^1(\mu^\w)= \min\{d_1,N\}
$$
for $\mu^\w$-almost all $x$.
\end{thm}

 \Proof Let $s$ be such that  (\ref{potest}) holds. By $(\ref{def_ua})$,  for all  $\bi \in I_\infty$ and  $0<r<1$,
\begin{eqnarray*}
\mu^{\w}(B(x^{\w}(\bi),r))    &=&\mu(\bj:\mid x^{\w}(\bi)-x^{\w}(\bj)\mid\leq r)                 \nonumber\\
&=&\mu\left(\bj:r^s |x^\w(\bi)-x^\w(\bj)|^{-s}
\geq 1\right)     \nonumber\\
&\leq& \int_{I_\infty}r^s |x^\w(\bi)-x^\w(\bj)|^{-s}
d\mu(\bj).
\end{eqnarray*}
Using (\ref{potest}),
\begin{eqnarray}
\E\big( \mu^{\w}(B(x^{\w}(\bi),r)) \big)
&\leq&\E \int_{I_\infty}r^s |x^\w(\bi)-x^\w(\bj)|^{-s}
d\mu(\bj)   
\nonumber \\
&\leq& r^s
\int_{I_\infty}\frac{c}{\phi^s(T_{\bi\wedge\bj})} d\mu(\bj)      \nonumber \\
&\leq& c r^s \sum_{k=1}^{\infty}\phi^s(T_{\bi\mid k})^{-1}
\mu(C_{\bi\mid k}). \label{series}
\end{eqnarray}
Since $\mu$ is an ergodic measure on $I_\infty$ and 
$\{\log(\phi^s(T_{\bi})^{-1} \mu(C_{\bi})\}_{\bi \in I}$ is superadditive,  the subadditive ergodic theorem(see \cite{Ste}) implies that $\frac{1}{k}\log\left(\phi^s(T_{\bi\mid k})^{-1} \mu(C_{\bi\mid k}) \right)\to  P(s,1)$ (given by (\ref{P1def})) for 
$\mu$-almost all $\bi \in I_\infty$.
Hence, 
if $s<d_1$,  then for $\mu$-almost all $\bi \in I_\infty$, we have 
$\lim_{k \to \infty}\left(\phi^s(T_{\bi\mid k})^{-1} \mu(C_{\bi\mid k})\right)^{1/k} \to 
\exp(P(s,1)) <1$, so the series (\ref{series}) converges. Thus if $t<s<d_1$
$$
\E \left( r^{-t} \mu^{\w}(B(x^{\w}(\bi),r))  \right)    \leq c_0 r^{s-t}
$$
for all $0<r<1$, for some $c_0$, for $\mu$-almost all $\bi$. 
Taking $r=2^{-l}$ and summing, 
\begin{eqnarray*}
\E \sum_{l=1}^\infty 2^{lt} \mu^{\w}(B(x^{\w}(\bi),2^{-l}))    \leq c_0 \sum_{l=1}^\infty  2^{-l(s-t)} < \infty.
\end{eqnarray*}
Since the limiting behaviour of the measure of balls is determined by the discrete radii $2^{-l}$, 
it follows that, for $\mu$-almost all $\bi \in I_\infty$, $$
r^{-t} \mu^{\w}(B(x^{\w}(\bi),r)) \to 0$$
as $r \to 0$, for almost all $\w$.
 Hence, for almost all
$\w$, 
\be
\liminf_{r\to 0 }\frac{\log \mu^{\w} (B(x,r))}{\log r}\geq t, \label{ineqt}
\ee
for $\mu^{\w}$-almost all $x$, and this is true for
all $t< \min\{d_1,N\}$, by choosing $s$ with  $t<s<$ for which (\ref{potest}) holds. Thus (\ref{ineqt}) holds with $t= \min\{d_1,N\}$, and the opposite estimate required comes from Proposition \ref{thm_ub}.

Finally, if  a probability measure $\nu$ is exact dimensional of dimension $d$ then $D^1(\nu)= d$; this  follows from the definitions (\ref{infdimint}) using the dominated convergence theorem, or see \cite{LSY}. Thus, in our case, $D^1(\mu^\w)=  \min\{d_1,N\}$ almost surely.
 \Endproof

%66666666666666666
\section{Specific cases}
\setcounter{equation}{0}
\setcounter{pro}{0}

We specialize Theorem \ref{thm_lb} to two cases of particular interest. First we consider a self-affine measure supported by a self-affine set. 
\begin{thm}\label{thmsa}
For $S_{i}(x) = T_{i}(x) +\w_{i}$ $(i=1,2,\ldots,m)$  let
$E^\w$ be the self-affine subset of $\bbbr^N$ satisfying
$(\ref{attractor})$. Let  $\mu$ be the Bernoulli measure on $I_\infty$ given by $(\ref{mP})$ and 
let $\mu^{\w}$ be the projection of $\mu$ onto $E^\w$, that is the self-affine measure satisfying $(\ref{msa})$.
Assume that $\|T_i\|<\frac{1}{2}$ for all $1\leq i \leq m$. Then,
for $Nm$-Lebesgue almost all $(\w_1, \ldots, \w_m) \in \bbbr^{Nm}$, the measure $\mu^\w$ is exact dimensional, with
\be
\lim_{r\to 0}\frac{\log \mu^\w (B(x,r))}{\log r} = D^1(\mu^\w)= \min\{d_1,N\}
\label{id1}
\ee
for $\mu^{\w}$-almost all $x$.
\end{thm}
\Proof 
It was shown in \cite[Lemma 2.2]{Fa1} for  $\|T_i\|<\frac{1}{3}$ and \cite{Sol} for  $\|T_i\|<\frac{1}{2}$ that, for all $\rho>0$ and $0<s<N$  with $s$ non-integral,
$$
\int_{\w_1,\ldots,\w_m\in B(0,\rho)} \frac{d\w_1\ldots d\w_m}{|x^\w(\bi)-x^\w(\bj)|^s}\leq
\frac{c}{\phi^s(T_{\bi\wedge \bj})},
$$
for all distinct $\bi$, $\bj\in I_\infty$, where $x^\w(\bi)$ is the point obtained by taking $ \w_ {i_1,\ldots,i_k} =  \w_{i_k}$ for all $(i_1,\ldots,i_k)$ in (\ref{points1})-(\ref{points3}). By taking $B$ a sufficiently large ball in $\bbbr^N$ we may ensure that (\ref{inclusions}) holds for all $\w_1,\ldots,\w_m\in B(0,\rho)$ for any given $\rho$. The conclusion follows from
Theorem \ref{thm_lb}, regarding normalised  $Nm$-dimensional Lebesgue measure on 
 $B(0,\rho)^m$ as a probability measure.
 \Endproof
 
 Note that in the setting of Theorem \ref{thmsa} we have almost surely that $D^q(\mu^\w)= \min\{d_q,N\}$ for $1<q \leq 2$, see \cite{Fa5}.

We now consider measures on statistically self-affine sets of the form considered in \cite{JPS,Fa6}. In particular this enables us to remove the restriction of $\|T_i\|<\frac{1}{2}$ of the previous theorem.

\begin{thm}\label{thmsaran}
Let
$E^\w$ be the almost self-affine subset of $\bbbr^N$ satisfying
$(\ref{construct})$. Let  $\mu$ be the Bernoulli measure on $I_\infty$ given by $(\ref{mP})$ and 
let $\mu^{\w}$ be the almost self-affine measure on $E^\w$ given by $(\ref{def_ua})$.
Suppose that $\|T_i\|<1$ for all $1\leq i \leq m$. Let $D$ be a bounded region in $ \bbbr^N$ and let $\P$ be a probability measure on $\Omega$ such that $\{\w_\bi: \bi \in \Omega\}$ are independent identically distributed random vectors in $D$ with a distribution that is absolutely continuous with
respect to $N$-dimensional Lebesgue measure. Then,
for $\P$-almost all $\w \in \Omega$, the measure $\mu^\w$ is exact dimensional with
\be
\lim_{r\to 0}\frac{\log \mu^\w (B(x,r))}{\log r} = D^1(\mu^\w)= \min\{d_1,N\}
\label{id2}
\ee
for $\mu^{\w}$-almost all $x$.
\end{thm}
\Proof
By taking $B$ a sufficiently large ball in $\bbbr^N$  we may ensure that (\ref{inclusions}) holds for all $\w$.
It was shown in \cite{JPS}, see also  \cite{Fa6}, that for this model, for all $0<s<N$  with $s$ non-integral,
$$\E\left( |x^\w(\bi)-x^\w(\bj)|^{-s}\right)\leq
\frac{c}{\phi^s(T_{\bi\wedge \bj})} 
$$
for all distinct $\bi$, $\bj\in I_\infty$, where $x^\w(\bi)$ is the random point given by (\ref{points1})-(\ref{points3}) for each $\bi \in I_\infty$.  The conclusion is immediate  from
Theorem \ref{thm_lb}.
 \Endproof
 
  Note that in the setting of Theorem \ref{thmsaran} we have that $D^q(\mu^\w)= \min\{d_q,N\}$ almost surely for all  $q>1$, see \cite{Fa6}.

%777777777777777
\section{Gibbs measures}
\setcounter{equation}{0}
\setcounter{pro}{0}

We may also get estimates for the local dimensions when   $\mu$ is an invariant Gibbs measure on $I_\infty$. Recall that a probability measure  $\mu$  on $I_\infty$ is a {\it Gibbs measure} if  there is a continuous $f: I_\infty \to \bbbr$ and  a real number $P(f)$, the {\it pressure} of $f$, such that for some $a>0$,
\be
a^{-1} \leq \frac{\mu(C_{\bi|k})}{\exp\big(-k P(f) + \sum_{j=0}^{k-1} f(\sigma^j (\bi)\big)}
\leq a \qquad (\bi \in I_\infty, k \in \mathbb{Z}^+),\label{gibbs}
\ee
where $\sigma$ is the shift on $I_\infty$ given by $\sigma(i_1,i_2,\ldots) = (i_2,i_3,\ldots)$. According to the variational principle, given $f: I_\infty \to \bbbr$  satisfying an 
$\epsilon$-H\"{o}lder condition, that is 
$|f(\bi) - f(\bj)|\leq  c d(\bi,\bj)^\epsilon$   for some $\epsilon>0, c>0$,  there exists an invariant Gibbs measure satisfying (\ref{gibbs}), and this provides a wide range of Gibbs measures. Note that  from (\ref{gibbs}),
\be
b^{-1}\mu(C_\bi)\mu(C_\bj)\leq \mu(C_{\bi\bj})\leq b\mu(C_\bi)\mu(C_\bj)
\qquad (\bi, \bj \in I), \label{quasimult}
\ee
where $b=a^3>0$. In particular $\{b\mu(C_\bi)\}_\bi$ is submultiplicative, and $\{b^{-1}\mu(C_\bi)\}_\bi$ is supermultiplicative, and this is enough to guarantee the existence of the limits in the definitions (\ref{Pdef}) and (\ref{P1def}), since, for example in the case of  $0\leq q<1$,  the limit 
$\lim_{k\to \infty}\frac{1}{k}
\log \sum_{|\bi|=k} \phi^s(T_\bi)^{1-q}b^q\mu(C_\bi)^q/(q-1)$ exists with a value unaltered if the $b^q$ term is dropped. As before, $P(s,q)$ is continuous in $s$ and strictly monotonic in $s$, so the definition (\ref{dqdef}) of $d_q \,(q \geq 0)$ remains valid in this case.

It follows from (\ref{gibbs})  that there are  constants $0< c_{-} \leq c_{+} <1$ such that 
$c_{-}^k \leq  \mu(C_\bi) \leq c_{+}^k$ for all $\bi\in I_k$. Thus Lemma \ref{lemPcts} remains true (taking $\beta_{-} = \log(\alpha_{+}^{-s}c_{-})$ and    
$\beta_{+} = \log(\alpha_{-}^{-s}c_{+})$ in (\ref{pcts})) except that we can no longer guarantee that $P(s,q)$ is lower semicontinuous at $q=1$.

Thus we get the following variants of Proposition \ref{thm_ub} and Theorem \ref{thm_lb}.

\begin{pro}\label{thm_ubg}
Let $\mu$  be a Gibbs measure on $I_\infty$. For any   $\w$ let $\mu^{\w}$ be the projection of $\mu$  onto $E^\w$ given by $(\ref{def_ua})$. Then
$$
\limsup_{r\to 0}\frac{\log\mu^\w(B(x,r))}{\log r}\leq \min\{\lim_{q\nearrow 1}d_q,N\} 
$$
for $\mu^\w$-almost all $x\in \bbbr^{N}$.
\end{pro}

\begin{thm}\label{thm_lbg}
Let $\mu$  be an invariant Gibbs measure on $I_\infty$ and
let $\mu^{\w}$ be the projection of $\mu$ onto $E^\w$. 
Suppose that there are numbers $s<\min\{d_1,N\}$ that are arbitrarily close to $\min\{d_1,N\}$ for which
there exists $c<\infty$ such that
$$
\E\left( |x^\w(\bi)-x^\w(\bj)|^{-s}\right)\leq
\frac{c}{\phi^s(T_{\bi\wedge \bj})} \quad (\bi \neq \bj \in I_\infty). 
$$
Then, for almost all $\w$,
\be
\min\{d_1,N\}\leq \lim_{r\to 0}\frac{\log \mu^\w (B(x,r))}{\log r}\leq 
 \min\{\lim_{q \nearrow 1} d_q,N\} \label{gibbsests}
\ee
for $\mu^\w$-almost all $x$.
\end{thm}

\Proof
The proof is similar to that of Theorem \ref{thm_lb}. For an invariant Gibbs measure, 
 $\{\log(\phi^s(T_{\bi})^{-1} b^{-1}\mu(C_{\bi})\}_{\bi \in I}$ is superadditive by (\ref{quasimult}), and since every invariant Gibbs measure is ergodic, the subadditive ergodic theorem gives that
$$\lim_{k \to \infty} \frac{1}{k}\log\left(\phi^s(T_{\bi\mid k})^{-1} \mu(C_{\bi\mid k}) \right)
=\lim_{k \to \infty} \frac{1}{k}\log\left(b^{-1}\phi^s(T_{\bi\mid k})^{-1} \mu(C_{\bi\mid k}) \right)=  P(s,1)$$
 for 
$\mu$-almost all $\bi \in I_\infty$.
Hence, 
if $s<d_1$,  then for $\mu$-almost all $\bi \in I_\infty$, we have that
$\lim_{k \to \infty}\left(\phi^s(T_{\bi\mid k})^{-1} \mu(C_{\bi\mid k})\right)^{1/k} \to 
\exp(P(s,1)) <1$, so the series (\ref{series}) converges. The proof concludes as in 
Theorem \ref{thm_lb}.
\Endproof

Of course, Theorem \ref{thm_lbg} may be specialized to self-similar measures and almost self-similars measures; thus if $\mu$ is a Gibbs measure rather than a Bernoulli measure in Theorems \ref{thmsa} and \ref{thmsaran}, the conclusions (\ref{id1}) and 
 (\ref{id2}) are replaced by (\ref{gibbsests}).
 
 For Gibbs measures it is not clear that  $P(s,q)$ need be lower semicontinuous in $q$ at $q=1$, so we cannot proceed as before to  get equality throughout (\ref{gibbsests}). It would be of interest to know when we do get equality here.

\end{document}